\documentclass[12pt]{article}
\usepackage[centertags]{amsmath}
\usepackage{amsfonts}
\usepackage{amssymb}
\usepackage{amsthm}
\usepackage{newlfont}
\usepackage{graphicx}

\newfont{\bb}{msbm10 at 12pt}

\newtheorem{theorem}{Theorem}[section]

\newtheorem{corollary}[theorem]{Corollary}

\begin{document}

\title{Some Characterizations of Rectifying Spacelike Curves in the Minkowski Space-Time}

\author{ Ahmad T. Ali\footnote{E-mail adress: atali71@yahoo.om (A. T. Ali).}\\Department of Mathematics\\
 Faculty of Science, Al-Azhar University\\
 Nasr City, 11448, Cairo, EGYPT\\
email: atali71@yahoo.com\\
\vspace*{1cm}\\
Mehmet \"{O}nder\\
Department of Mathematics\\
Faculty of Science and Arts, Celal Bayar University\\
Muradiye Campus, Manisa, TURKEY\\
email: mehmet.onder@bayar.edu.tr}
\date{}

\maketitle
\begin{abstract} In this paper, we define a rectifying spacelike curve in the Minkowski space-time $E_1^4$ as a curve whose position vector always lies in orthogonal complement $N^{\bot}$ of its principal normal vector field $N$. In particular, we study the rectifying spacelike curves in $E_1^4$ and characterize such curves in terms of their curvature functions.
\end{abstract}



\section{Introduction}
In the Euclidean space $E^3$, rectifying curves are introduced by Chen in \cite{chen1}
as space curves whose position vector always lies in its rectifying plane, spanned by the tangent and the binormal vector fields $\vec{T}$ and $\vec{B}$ of the curve. Therefore, the position vector $\vec{\alpha}$ of a rectifying curve satisfies the equation
$$
\vec{\alpha}(s)=\lambda(s)\vec{T}(s)+\mu(s)\vec{B}(s),
$$
for some differentiable functions $\lambda$ and $\mu$ in arclength function $s$. The Euclidean rectifying curves are studied in \cite{chen1, chen2}. In particular, it is shown in \cite{chen2} that there exists a simple relationship between the rectifying curves and the centrodes, which play some important roles in mechanics, kinematics as well as in differential geometry in defining the curves of constant precession. The rectifying curves are also studied in \cite{chen2} as the extremal curves. In the Minkowski 3-space $E_1^3$, the rectifying curves are investigated in \cite{ilar2}. The rectifying curves are also studied in \cite{ilar3} as the centrodes and extremal curves. In the Euclidean 4-space $E^4$, the rectifying curves are investigated in \cite{ilar1}.

In analogue with the rectifying curve the curve whose position vector always lies in its normal plane spanned by the principal normal and the binormal vector fields $\vec{N}$ and $\vec{B}$ of the curve is called normal curve in Euclidean 3-space $E^3$ and it is well known that normal curves are spherical curves in $E^3$[2]. Similar definition and characterizations of spacelike, timelike (and also null) and dual timelike normal curves are given in references [7], [8] and [10], respectively. The spacelike normal curve in the Minkowski 4-space $E_1^4$ can be defined as a curve whose position vector always lies in the orthogonal complement $\vec{T}^{\bot}$ of its tangent vector field $\vec{T}$ which is given by
$$
\vec{T}^{\bot}=\{\vec{W}\in E_1^4\mid g(\vec{W},\vec{T})=0\}.
$$
In [1], Camci and others have shown that a spacelike curve is a hyperbolic spherical curve iff the following equation holds
$$
\vec{\alpha}-m=-(1/k_1)\vec{N}-(1/k_2)(1/k_1)'\vec{B_1}+(1/k_3)[k_2/k_1+((1/k_2)(1/k_1)')']\vec{B_2},
$$
where $m$ is constant, $k_1$, $k_2$ and $k_3$ are the first, the second and the third curvatures of the curve $\alpha$, respectively. By using the definition of spacelike normal curves in Minkowski 4-space $E_1^4$ and the last equality, we can say that spacelike normal curves are hyperbolic spherical curves in Minkowski 4-space $E_1^4$.\\

In this paper, in analogy with the Minkowski 3-dimensional case, we define the rectifying curve in the Minkowski 4-space $E_1^4$ as a curve whose position vector always lies in the orthogonal complement $\vec{N}^{\bot}$ of its principal normal vector field $\vec{N}$. Consequently, $\vec{N}^{\bot}$ is given by
$$
\vec{N}^{\bot}=\{\vec{W}\in E_1^4\mid g(\vec{W},\vec{N})=0\},
$$
where $g(.,.)$ denotes the standard pseudo scalar product in $E_1^4$. Hence $\vec{N}^{\bot}$ is a 3-dimensional subspace of $E_1^4$, spanned by the tangent, the first binormal and the second binormal vector fields $\vec{T}$, $\vec{B}_1$ and $\vec{B}_2$ respectively. Therefore, the position vector with respect to some chosen origin, of a rectifying spacelike curve $\vec{\alpha}$ in Minkowski space-time $E_1^4$, satisfies the equation
\begin{equation}\label{equi1}
\vec{\alpha}(s)=\lambda(s)\vec{T}(s)+\mu(s)\vec{B}_1(s)+\nu(s)\vec{B}_2,
\end{equation}
for some differentiable functions $\lambda(s)$, $\mu(s)$ and $\nu(s)$ in arclength function $s$. Next, characterize rectifying spacelike curves in terms of their curvature functions $\kappa_1(s)$, $\kappa_2(s)$ and $\kappa_3(s)$ and give the necessary and the sufficient conditions for arbitrary curve in $E_1^4$ to be a rectifying. Moreover, we obtain an explicit equation of a rectifying spacelike curve in $E_1^4$ and give the relation between rectifying and normal spacelike curves in $E_1^4$.

\section{Preliminaries}

In this section, we prepare basic notations on Minkowski space-time $E_1^4$. Let $\vec{\alpha}:I\subset R\rightarrow E_1^4$ be arbitrary curve in the Minkowski space-time $E_1^4$. Recall that the curve $\vec{\alpha}$ is said to be unit speed(or parameterized by arclength function $s$) if $g(\vec{\alpha}^{\prime},\vec{\alpha}^{\prime})=\pm1$, where $g(.,.)$ denotes the standard pseudo scalar product in $E_1^4$ given by
$$
g(\vec{v},\vec{w})=-v_1w_1+v_2w_2+v_3w_3+v_4w_4,
$$
for each $\vec{v}=(v_1,v_2,v_3,v_4)\in E_1^4$ and $\vec{w}=(w_1,w_2,w_3,w_4)\in E_1^4$. An arbitrary vector $\vec{v}\in E_1^4$ can have one of three Lorentzian causal characters; it can be spacelike if $g(\vec{v},\vec{v})>0$ or $\vec{v}=0$, timelike if $g(\vec{v},\vec{v})<0$ and null (lightlike) if $g(\vec{v},\vec{v})=0$ and $\vec{v}\neq0$. Similarly, an arbitrary curve $\vec{\alpha}=\vec{\alpha}(s)$ can locally be spacelike, timelike or null (lightlike), if all of its velocity vectors $\vec{\alpha}^{\prime}(s)$ are respectively spacelike, timelike or null (lightlike). Also recall that the pseudo-norm of an arbitrary vector $\vec{v}\in E_1^4$ is given by $\parallel\vec{v}\parallel=\sqrt{\mid g(\vec{v},\vec{v})\mid}$. Therefore $\vec{v}$ is a unit vector if $g(\vec{v},\vec{v})=\pm1$. The velocity of the curve $\vec{\alpha}(s)$ is given by $\parallel\vec{\alpha}^{\prime}(s)\parallel$. Next, vectors $\vec{v},\vec{w}$ in $E_1^4$ are said to be orthogonal if $g(\vec{v},\vec{w})=0$.

Denote by $\{\vec{T}(s),\vec{N}(s),\vec{B}_1(s),\vec{B}_2(s)\}$ the moving Frenet frame along the curve $\vec{\alpha}(s)$ in the space $E_1^4$, where $\vec{T}(s)$, $\vec{N}(s)$, $\vec{B}_1(s)$ and $\vec{B}_2(s)$ are the tangent, principal normal, the first binormal and second binormal fields, respectively. For an arbitrary spacelike curve $\vec{\alpha}(s)$ with spacelike principal normal $\vec{N}$ in the space $E_1^4$, the following Frenet formula are given in \cite{walr1, petr1, camc1, kaza1, turg1, yilm1}:
\begin{equation}\label{equi2}
 \left[
   \begin{array}{c}
     \vec{T}^{\prime} \\
     \vec{N}^{\prime} \\
     \vec{B}_1^{\prime} \\
     \vec{B}_2^{\prime}\\
   \end{array}
 \right]=\left[
           \begin{array}{cccc}
             0 & \kappa_1 & 0 & 0 \\
             -\kappa_1 & 0 & \kappa_2 & 0 \\
             0 & -\varepsilon\kappa_2 & 0 & \kappa_3 \\
             0 & 0 & \kappa_3 & 0 \\
           \end{array}
         \right]\left[
   \begin{array}{c}
     \vec{T} \\
     \vec{N} \\
     \vec{B}_1 \\
     \vec{B}_2\\
   \end{array}
 \right],
 \end{equation}
where
\begin{equation}\label{equi99}
g(\vec{B}_1,\vec{B}_1)=-g(\vec{B}_2,\vec{B}_2)=\varepsilon=\pm1,\,\,\,\,\,
g(\vec{T},\vec{T})=g(\vec{N},\vec{N})=1.
\end{equation}
Recall the functions $\kappa_1(s)$, $\kappa_2(s)$ and $\kappa_3(s)$ are called respectively, the first, the second and the third curvatures of timelike curve $\vec{\alpha}(s)$. Here, $\varepsilon$ determines the kind of spacelike curve $\alpha(s)$. If $\varepsilon=1$, then $\alpha(s)$ is a spacelike curve with spacelike first binormal $\vec{B}_1$ and timelike second binormal $B_2$. If $\varepsilon=-1$, then $\alpha(s)$ is a spacelike curve with timelike first binormal $\vec{B}_1$ and spacelike second binormal $B_2$. If $\kappa_3(s)\neq0$ for each $s\in I\subset R$, the curve $\vec{\alpha}$ lies fully in $E_1^4$. Recall that the hyperbolic unit sphere $H_0^3(1)$ in $E_1^4$, centered at the origin, is the hypersurface defined by
$$
H_0^3(1)=\{\vec{X}\in E_1^4\mid g(\vec{X},\vec{X})=-1\}.
$$

\section{Some characterizations of rectifying curves in $E_1^4$}

In this section, we firstly characterize the rectifying spacelike curves in Minkowski space-time in terms of their curvatures. Let $\vec{\alpha}=\vec{\alpha}(s)$ be a unit speed rectifying spacelike curve in $E_1^4$, with non-zero curvatures $\kappa_1(s)$, $\kappa_2(s)$ and $\kappa_3(s)$. By definition, the position vector of the curve $\vec{\alpha}$ satisfies the equation (\ref{equi1}), for some differentiable functions $\lambda(s)$, $\mu(s)$ and $\nu(s)$. Differentiating the equation (\ref{equi1}) with respect to $s$ and using the Frenet equations (\ref{equi2}), we obtain
$$
\vec{T}=\lambda^{\prime}\vec{T}+(\lambda\kappa_1-\varepsilon\mu\kappa_2)\vec{N}+(\mu^{\prime}+\nu\kappa_3)\vec{B}_1
+(\nu^{\prime}+\mu\kappa_3)\vec{B}_2.
$$
It follows that
\begin{equation}\label{equi4}
\begin{array}{l l}
\lambda^{\prime}=1,\\
\lambda\kappa_1-\varepsilon\mu\kappa_2=0,\\
\mu^{\prime}+\nu\kappa_3=0,\\
\nu^{\prime}+\mu\kappa_3=0,
\end{array}
\end{equation}
and therefore
\begin{equation}\label{equi5}
\begin{array}{l l}
\lambda=s+c,\\
\mu=\varepsilon\,\dfrac{\kappa_1(s)(s+c)}{\kappa_2},\\
\nu=-\varepsilon\,\dfrac{\kappa_1(s)\kappa_2(s)+(s+c)(\kappa_1^{\prime}(s)\kappa_2(s)
-\kappa_1(s)\kappa_2^{\prime}(s))}{\kappa_2^2(s)\kappa_3(s)},
\end{array}
\end{equation}
where $c\in R$. In this way the functions $\lambda(s)$, $\mu(s)$ and $\nu(s)$ are expressed in terms of the curvature functions $\kappa_1(s)$, $\kappa_2(s)$ and $\kappa_3(s)$ of the curve $\alpha(s)$. Moreover, by using the last equation in (\ref{equi4}) and relation (\ref{equi5}) we easily find that the curvatures $\kappa_1(s)$, $\kappa_2(s)$ and $\kappa_3(s)$satisfy the equation
\begin{equation}\label{equi6}
\begin{array}{l l}
\dfrac{\kappa_1(s)\kappa_3(s)(s+c)}{\kappa_2(s)}-
\Big[\dfrac{\kappa_1(s)\kappa_2(s)+(s+c)(\kappa_1^{\prime}(s)\kappa_2(s)
-\kappa_1(s)\kappa_2^{\prime}(s))}{\kappa_2^2(s)\kappa_3(s)}\Big]^{\prime}=0,\,\,\,c\in R.
\end{array}
\end{equation}
The condition (\ref{equi6}) can be written as:
\begin{equation}\label{equi8}
\begin{array}{l l}
\dfrac{\kappa_1(s)(s+c)}{\kappa_2(s)}-
\dfrac{1}{\kappa_3(s)}\dfrac{d}{ds}\Big[\dfrac{1}{\kappa_3(s)}\dfrac{d}{ds}\Big(\dfrac{\kappa_1(s)(s+c)
}{\kappa_2(s)}\Big)\Big]=0.
\end{array}
\end{equation}
If we change the variable $s$ by the variable $t$ as the following
$$
\dfrac{d}{dt}=\dfrac{1}{\kappa_3(s)}\dfrac{d}{ds}\,\,\,\,\text{or}\,\,\,\,t=\int_0^s\kappa_3(s)ds,
$$
the equation (\ref{equi8}) takes the following form
\begin{equation}\label{equi9}
\begin{array}{l l}
\dfrac{\kappa_1(s)(s+c)}{\kappa_2(s)}-
\dfrac{d^2}{dt^2}\Big[\dfrac{\kappa_1(s)(s+c)
}{\kappa_2(s)}\Big]=0.
\end{array}
\end{equation}
General solution of this equation is
\begin{equation}\label{equi10}
\begin{array}{l l}
\dfrac{\kappa_1(s)(s+c)}{\kappa_2(s)}=\varepsilon\Big(A\cosh\int_0^s\kappa_3(s)ds+B\sinh\int_0^s\kappa_3(s)ds\Big),
\end{array}
\end{equation}
where $A$ and $B$ are arbitrary constants. Then from (\ref{equi5}) we have
\begin{equation}\label{equi11}
\begin{array}{l l}
\lambda(s)=s+c\\
\mu(s)=A\cosh\int_0^s\kappa_3(s)ds+B\sinh\int_0^s\kappa_3(s)ds\\
\nu(s)=-\Big(A\sinh\int_0^s\kappa_3(s)ds+B\cosh\int_0^s\kappa_3(s)ds\Big).
\end{array}
\end{equation}

Conversely, assume that the curvatures $\kappa_1(s)$, $\kappa_2(s)$ and $\kappa_3(s)$, of an arbitrary unit speed spacelike curve in $E_1^4$, satisfy the equation (\ref{equi10}). Let us consider the vector $\vec{X}\in E_1^4$ given by
\begin{equation}\label{equi7}
\begin{array}{l l}
\vec{X}(s)&=\vec{\alpha}(s)-(s+c)\vec{T}(s)-\Big(A\cosh\int_0^s\kappa_3(s)ds+B\sinh\int_0^s\kappa_3(s)ds\Big)\vec{B}_1(s)\\
&+\Big(A\sinh\int_0^s\kappa_3(s)ds+B\cosh\int_0^s\kappa_3(s)ds\Big)\vec{B}_2(s).
\end{array}
\end{equation}
By using the relation (\ref{equi2}) and (\ref{equi10}), we easily find $\vec{X}^{\prime}(s)=0$, which means that $\vec{X}$ is a constant vector. This implies that $\alpha(s)$ is congruent to a rectifying curve. In this way, the following theorem is proved.

\begin{theorem}\label{th-31} Let $\vec{\alpha}(s)$ be unit speed spacelike curve in $E_1^4$, with non-zero curvatures $\kappa_1(s)$, $\kappa_2(s)$ and $\kappa_3(s)$. Then $\vec{\alpha}(s)$ is congruent to a rectifying spacelike curve if and only if
$$
\dfrac{\kappa_1(s)(s+c)}{\kappa_2(s)}=\varepsilon\Big(A\cosh\int_0^s\kappa_3(s)ds+B\sinh\int_0^s\kappa_3(s)ds\Big).
$$
\end{theorem}

In particular, assume that all curvature functions $\kappa_1(s)$, $\kappa_2(s)$ and $\kappa_3(s)$ of rectifying spacelike curve $\vec{\alpha}$ in $E_1^4$, are constant and different from zero. Then equation (\ref{equi6}) easily implies a contradiction. Hence we obtain the following theorem.

\begin{theorem}\label{th-32} There are no rectifying spacelike curves lying in $E_1^4$, with non-zero constant curvatures $\kappa_1(s)$, $\kappa_2(s)$ and $\kappa_3(s)$.
\end{theorem}

In the next theorem, we give the necessary and the sufficient conditions for the spacelike curve $\alpha(s)$ in $E_1^4$ to be a rectifying curve.

\begin{theorem}\label{th-33} Let $\alpha(s)$ be unit speed rectifying spacelike curve in $E_1^4$, with non-zero curvatures $\kappa_1(s)$, $\kappa_2(s)$ and $\kappa_3(s)$. Then the following statements hold:

{\bf{(i)}} The distance function $\rho(s)=\parallel\vec{\alpha}(s)\parallel$ satisfies $\rho^2(s)=s^2+c_1s+c_2$, $c_1\in R$ and $c_2\in R$.

{\bf{(ii)}} The tangential component of the position vector of the rectifying spacelike curve is given by $g(\vec{\alpha}(s),\vec{T}(s))=s+c$, $c\in R$.

{\bf{(iii)}} The normal component $\vec{\alpha}^{N}(s)$ of the position vector of the rectifying spacelike curve has constant length and the distance function $\rho(s)$is non-constant.

{\bf{(iv)}} The first binormal component and the second binormal component of the position vector of the rectifying spacelike curve are respectively given by
\begin{equation}\label{equi21}
\begin{array}{l l}
g(\vec{\alpha}(s),\vec{B}_1(s))=\varepsilon\Big(A\cosh\int_0^s\kappa_3(s)ds+B\sinh\int_0^s\kappa_3(s)ds\Big)\\
g(\vec{\alpha}(s),\vec{B}_2(s))=\varepsilon\Big(A\sinh\int_0^s\kappa_3(s)ds-B\cosh\int_0^s\kappa_3(s)ds\Big).
\end{array}
\end{equation}

Conversely, if $\vec{\alpha}(s)$ is a unit speed curve in $E_1^4$ with non-zero curvatures $\kappa_1(s)$, $\kappa_2(s)$ and $\kappa_3(s)$ and one of the statements (i), (ii), (iii) or (iv) holds, then $\vec{\alpha}(s)$ is a rectifying spacelike curve.
\end{theorem}

{\bf{Proof.}} Let us first suppose that $\vec{\alpha}(s)$ is a unit speed rectifying spacelike curve in $E_1^4$ with non-zero curvatures $\kappa_1(s)$, $\kappa_2(s)$ and $\kappa_3(s)$. The position vector of the curve $\vec{\alpha}(s)$ satisfies the equation (\ref{equi1}), where the functions $\lambda(s)$, $\mu(s)$ and $\nu(s)$ satisfy relation (\ref{equi11}). From relation (\ref{equi1}) and (\ref{equi11}) we have
\begin{equation}
\begin{array}{l l}
g(\vec{\alpha},\vec{\alpha})&=\lambda^2+\varepsilon\Big(\mu^2(s)-\nu^2(s)\Big),\\
&=(s+c)^2+\varepsilon(A^2-B^2).
\end{array}\nonumber
\end{equation}
Therefore, $\rho^2(s)=s^2+c_1s+c_2$, $c_1\in R$ and $c_2\in R$, which proves statement (i).

But using the relations (\ref{equi1}) and (\ref{equi5}) we easily get $g(\vec{\alpha}(s),\vec{T}(s))=s+c$, $c\in R$, so the statement (ii) is proved.

Note that the position vector of an arbitrary curve $\vec{\alpha}(s)$ in $E_1^4$ can be decomposed as $\vec{\alpha}(s)=m(s)\vec{T}(s)+\vec{\alpha}^{N}(s)$, where $m(s)$ is arbitrary differentiable function and $\vec{\alpha}^{N}(s)$ is the normal component of the position vector. If $\vec{\alpha}(s)$ is a rectifying spacelike curve, relation (\ref{equi1}) implies $\vec{\alpha}^{N}(s)=\mu(s)\vec{B}_1(s)+\nu(s)\vec{B}_2(s)$ and therefore $g(\vec{\alpha}^{N}(s),\vec{\alpha}^{N}(s))=\varepsilon\Big(\mu^2(s)-\nu^2(s)\Big)$. Moreover, by using (\ref{equi11}), we find $\parallel\vec{\alpha}^{N}(s)\parallel=\varepsilon(A^2-B^2)=a$, $a\in R$. By statement (i), $\rho(s)$ is non-constant function, which proves statement (iii).

Finally, using (\ref{equi1}), (\ref{equi99}) and (\ref{equi11}) we easily obtain (\ref{equi21}), which proves statement (iv).

Conversely, assume that statement (i) holds. Then $g(\vec{\alpha}(s),\vec{\alpha}(s))=s^2+c_1s+c_2$, $c_1\in R$, $c_2\in R$. By differentiating the previous equation two times with respect to $s$ and using (\ref{equi2}), we obtain $g(\vec{\alpha}(s),\vec{N}(s))=0$, which implies that $\vec{\alpha}$ is a rectifying spacelike curve.

If statement (ii) holds, in a similar way it follows that $\vec{\alpha}$ is a rectifying spacelike curve.

If statement (iii) holds, let us put $\vec{\alpha}(s)=m(s)\vec{T}(s)+\vec{\alpha}^{N}(s)$, where $m(s)$ is arbitrary differentiable function. Then
$$
g(\vec{\alpha}^{N}(s),\vec{\alpha}^{N}(s))=g(\vec{\alpha}(s),\vec{\alpha}(s))-2g(\vec{\alpha}(s),\vec{T}(s))m(s)+m^2(s).
$$
Since $g(\vec{\alpha}(s),\vec{T}(s))=m(s)$, it follows that
$$
g(\vec{\alpha}^{N}(s),\vec{\alpha}^{N}(s))=g(\vec{\alpha}(s),\vec{\alpha}(s))-g(\vec{\alpha}(s),\vec{T}(s))^2,
$$
where $g(\vec{\alpha}(s),\vec{\alpha}(s))=\rho^2(s)\neq$ constant. Defferentiating the previous equation with respect to $s$ and using (\ref{equi2}), we find
$$
\kappa_1(s)\,g(\vec{\alpha}(s),\vec{T}(s))\,g(\vec{\alpha}(s),\vec{N}(s))=0.
$$
It follows that $g(\vec{\alpha}(s),\vec{N}(s))=0$ and hence the spacelike curve $\vec{\alpha}$ is a rectifying.

If the statement (iv) holds, by taking the derivative of the equation
$$
g(\vec{\alpha}(s),\vec{B}_1(s))=\varepsilon\Big(A\cosh\int_0^s\kappa_3(s)ds+B\sinh\int_0^s\kappa_3(s)ds\Big),
$$
with respect to $s$ and using (\ref{equi2}), we obtain
$$
-\varepsilon\kappa_2(s)g(\vec{\alpha}(s),\vec{N}(s))+\kappa_3g(\vec{\alpha}(s),\vec{B}_2(s))
=\varepsilon\kappa_3\Big(A\sinh\int_0^s\kappa_3(s)ds+B\cosh\int_0^s\kappa_3(s)ds\Big).
$$
By using (\ref{equi21}), the last equation becomes $g(\vec{\alpha}(s),\vec{N}(s))=0$, which means that $\vec{\alpha}$ is a rectifying spacelike curve. This proves the theorem.

In the next theorem, we find the parametric equation of a rectifying curve.

\begin{theorem}\label{th-34} Let $\alpha:I\subset R\rightarrow E_1^4$ be a spacelike curve in $E_1^4$ given by $\vec{\alpha}(t)=\rho(t)\vec{y}(t)$, where $\rho(t)$ is arbitrary positive function and $\vec{y}(t)$ is a unit speed curve in the hyperbolic unit sphere $H_0^3(1)$. Then $\vec{\alpha}$ is a rectifying spacelike curve if and only if
\begin{equation}\label{equi22}
\begin{array}{l l}
\rho(t)=\dfrac{a}{\sin(t+t_0)},\,\,\,\,\, a\in R_0,\,\,\,\,\,t_0\in R.
\end{array}
\end{equation}
\end{theorem}

{\bf{Proof.}} Let $\vec{\alpha}$ be a curve in $E_1^4$ given by
$$
\vec{\alpha}(t)=\rho(t)\vec{y}(t)
$$
where $\rho(t)$ is arbitrary positive function and $\vec{y}(t)$ is a unit speed curve in the hyperbolic unit sphere $H_0^3(1)$. By taking the derivative of the previous equation with respect to $t$, we get
$$
\vec{\alpha}^{\prime}(t)=\rho^{\prime}(t)\vec{y}(t)+\rho(t)\vec{y}^{\,\prime}(t).
$$
Hence the unit tangent vector of $\vec{\alpha}$ is given by
\begin{equation}\label{equi23}
\begin{array}{l l}
\vec{T}=\dfrac{\rho^{\prime}(t)}{v(t)}\vec{y}(t)+\dfrac{\rho(t)}{v(t)}\vec{y}^{\,\prime}(t),
\end{array}
\end{equation}
where $v(t)=\parallel\vec{\alpha}^{\prime}(t)\parallel$ is the speed of $\vec{\alpha}$. Differentiating the equation (\ref{equi23}) with respect to $t$, we find
\begin{equation}\label{equi24}
\begin{array}{l l}
\vec{T}^{\prime}=\Big(\dfrac{\rho^{\prime}}{v}\Big)^{\prime}\vec{y}+\Big(\dfrac{2\rho^{\prime}}{v}-
\dfrac{\rho\rho^{\prime}(\rho+\rho^{\prime\prime})}{v^3}\Big)\vec{y}^{\,\prime}
+\Big(\dfrac{\rho}{v}\Big)\vec{y}^{\,\prime\prime}.
\end{array}
\end{equation}
Let $\vec{Y}$ be the unit vector field in $E_1^4$ satisfying the equations $g(\vec{Y},\vec{y})=g(\vec{Y},\vec{y}^{\,\prime})=g(\vec{Y},\vec{y}\times\vec{y}^{\,\prime})=0$. Then $\{\vec{y},\vec{y}^{\,\prime},\vec{y}\times\vec{y}^{\,\prime},\vec{Y}\}$ is orthogonal frame of $E_1^4$. Therefore decomposition of $\vec{y}^{\,\prime\prime}$ with respect the frame $\{\vec{y},\vec{y}^{\,\prime},\vec{y}\times\vec{y}^{\,\prime},\vec{Y}\}$ reads
\begin{equation}\label{equi25}
\begin{array}{l l}
\vec{y}^{\,\prime\prime}=g(\vec{y}^{\,\prime\prime},\vec{y})\,\vec{y}
+g(\vec{y}^{\,\prime\prime},\vec{y}^{\,\prime})\,\vec{y}^{\,\prime}
+g(\vec{y}^{\,\prime\prime},\vec{y}\times\vec{y}^{\,\prime})\,\vec{y}\times\vec{y}^{\,\prime}
+g(\vec{y}^{\,\prime\prime},\vec{Y})\,\vec{Y}.
\end{array}
\end{equation}
Since $g(\vec{y},\vec{y})=-1$ and $g(\vec{y}^{\,\prime},\vec{y}^{\,\prime})=1$, it follows that $g(\vec{y}^{\,\prime\prime},\vec{y})=-1$ and $g(\vec{y}^{\,\prime\prime},\vec{y}^{\,\prime})=0$, so the equation (\ref{equi25}) becomes
\begin{equation}\label{equi26}
\begin{array}{l l}
\vec{y}^{\,\prime\prime}=-\vec{y}
+g(\vec{y}^{\,\prime\prime},\vec{y}\times\vec{y}^{\,\prime})\,\vec{y}\times\vec{y}^{\,\prime}
+g(\vec{y}^{\,\prime\prime},\vec{Y})\,\vec{Y}.
\end{array}
\end{equation}
Substituting (\ref{equi25}) into (\ref{equi24}) and applying Frenet formulas for arbitrary speed curves in $E_1^4$, we find
\begin{equation}\label{equi27}
\begin{array}{l l}
\kappa_1v\vec{N}&=\Big[\Big(\dfrac{\rho^{\prime}}{v}\Big)^{\prime}-\dfrac{\rho}{v}\Big]\vec{y}
+\Big(\dfrac{2\rho^{\prime}}{v}-
\dfrac{\rho\rho^{\prime}(\rho+\rho^{\prime\prime})}{v^3}\Big)\vec{y}^{\,\prime}
+\dfrac{g(\vec{y}^{\,\prime\prime},\vec{y}\times\vec{y}^{\,\prime})}{v}\,\vec{\alpha}\times\vec{y}^{\,\prime}\\
&+\Big(\dfrac{\rho}{v}\Big)g(\vec{y}^{\,\prime\prime},\vec{Y})\vec{Y}.
\end{array}
\end{equation}
Since $g(\vec{y},\vec{y})=-1$, we have $g(\vec{y},\vec{y}^{\,\prime})=0$ and thus $g(\vec{\alpha},\vec{y}^{\,\prime})=0$. We also have $g(\vec{\alpha},\vec{Y})=0$. By definition, $\vec{\alpha}$ is a rectifying spacelike curve in $E_1^4$ if and only if $g(\vec{\alpha},\vec{N})=0$. Therefore, after taking the scalar product of (\ref{equi27}) with $\vec{\alpha}$, we have $g(\vec{\alpha},\vec{N})=0$ if and only if
\begin{equation}\label{equi28}
\begin{array}{l l}
\Big(\dfrac{\rho^{\prime}}{v}\Big)^{\prime}-\dfrac{\rho}{v}=0,
\end{array}
\end{equation}
whose general solutions are given by (\ref{equi22}). This proves the theorem.\\

In Theorem 3.4, since $\vec{y}(s)$ is a unit speed spacelike curve in the hyperbolic unit sphere $H_0^3(1)$, $\vec{y}(s)$ is a spacelike normal curve in Minkowski 4-space $E_1^4$. So, Theorem 3.4 give the relation between spacelike rectifying and spacelike normal curves in Minkowski 4-space $E_1^4$. Then we can give the following corollary.\\

\begin{corollary} \label{cr-31} In Minkowski 4-space $E_1^4$, the construction of the rectifying curves can be made by using the spacelike normal curves.
\end{corollary}

{\bf Example:} Let us consider the curve $\vec{\alpha}(s)=\dfrac{a}{\sin(s+s_0)}\Big(\cosh(s),0,\sinh(s),0\Big)$, $a\in R_0$, $s_0\in R$ in $E_1^4$. This curve has a form $\vec{\alpha}(s)=\rho(s)\vec{y}$, where $\rho(s)=\dfrac{a}{\,\sin(s+s_0)}$ and $\vec{y}(s)=\Big(\cosh(s),0,\sinh(s),0\Big)$. Since $g(\vec{y}(s),\vec{y}(s))=-1$ and $\parallel\vec{y}^{\,\prime}(s)\parallel=1$, $\vec{y}(s)$ is a unit speed spacelike curve in the hyperbolic unit sphere $H_0^3(1)$. According to theorem \ref{th-34}, $\vec{\alpha}$ is a rectifying spacelike curve lying fully in $E_1^4$.


\begin{thebibliography}{99}

\bibitem{camc1} Camci C, \.{I}larslan K, \v{S}u\'{c}urovi\'{c} E. On pseudohyperbolical curves in Minkowski space-time. Turk J Math 2003;27:315-328.

\bibitem{chen1} Chen BY. When does the position vector of a space curve always lie in its rectifying plane?.
Amer Math Monthly 2003;110:147-152.

\bibitem{chen2} Chen BY, Dillen F. Rectifying curves as centrodes and extremal curves.
Bull Inst Math Acadimia Sinica 2005;2:77-90.

\bibitem{ilar1} \.{I}larslan K, Ne\v{s}ovi\'{c} E. Some characterizations of rectifying curves in the Euclidean space $E^4$. Turk J Math 2008;32:21-30.

\bibitem{ilar2} \.{I}larslan K, Ne\v{s}ovi\'{c} E, Petrovi\'{e}-Torga\v{s}ev M. Some characterizations of rectifying curves in Minkowski 3-space. Novi Sad J Math 2003;33(2):23-32.

\bibitem{ilar3} \.{I}larslan K, Ne\v{s}ovi\'{c} E. On rectifying curves as centrodes and extremal curves in the Minkowski 3-space. Novi Sad J Math 2007;37(1):53-64.

\bibitem{ilar4} \.{I}larslan K. Spacelike Normal Curves in Minkowski space $E_1^3$. Turk J Math 2005;29:53-63.

\bibitem{ilar5} \.{I}larslan K, Ne\v{s}ovi\'{c} E. Timelike and null Normal Curves in Minkowski space. Indian J Pure and Appl. Math 2004;35(7):881-888.

\bibitem{kaza1} Kazaz M, \"{O}nder M, Kocayi\v{g}it H. Spacelike Curves of constant Breadth in Minkowski 4-space. Int J Math Anal 2008;2(22):1061-1068.

\bibitem{"{O}nde1} \"{O}nder M. Dual timelike normal and dual timelike spherical curves in dual Minkowski space $D_1^3$. SDU Journal of Science 2006;1 (1-2):77-86.

\bibitem{petr1} Petrovi\'{c}-Torga\v{s}ev M, \v{S}u\'{c}urovi\'{c} E. W-Curves in Minkowski space-time. Novi Sad J Math 2002;32(2):55-65.

\bibitem{turg1} Turgut M, Yilmaz S. On the Frenet frame and a characterization of space-like involute-evolute curve couple in Minkowski space-time. Int Math Forum 2008;3(16):793-801.

\bibitem{walr1} Walrave J. Curves and surfaces in Minkowski space, doctoral thesis, K U Leuven, Fac Sci, Leuven, 1995.

\bibitem{yilm1} Yilmaz S, Turgut M. Relations among Frenet Apparatus of space-like Bertrand W-curve couples in Minkowski space-time. Int Math Forum 2008;3(32):1575-1580.

\end{thebibliography}
\end{document}